\magnification1200
\input amssym.def 
\input amssym.tex 
\def\SetAuthorHead#1{}
\def\SetTitleHead#1{}
\def\NoindentAfter{\everypar={\setbox0=\lastbox\everypar={}}}
\def\H#1\par#2\par{{\baselineskip=15pt\parindent=0pt\parskip=0pt
 \leftskip= 0pt plus.2\hsize\rightskip=0pt plus.2\hsize
 \bf#1\unskip\break\vskip 4pt\rm#2\unskip\break\hrule
 \vskip40pt plus4pt minus4pt}\NoindentAfter}
\def\HH#1\par{{\bigbreak\noindent\bf#1\medbreak}\NoindentAfter}
\def\HHH#1\par{{\bigbreak\noindent\bf#1\unskip.\kern.4em}}
\def\th#1\par{\medbreak\noindent{\bf#1\unskip.\kern.4em}\it}
\def\endth{\medbreak\rm}
\def\pf#1\par{\medbreak\noindent{\it#1\unskip.\kern.4em}}
\def\df#1\par{\medbreak\noindent{\it#1\unskip.\kern.4em}}

\let\Roster\bgroup\let\endRoster\egroup
\def\\{}\def\text#1{\hbox{\rm #1}}

\def\MaxReferenceTag#1{}
\def\qedbox{\vrule width2mm height2mm\hglue1mm\relax}
\def\qed{\ifmmode\qedbox\else\hglue5mm\unskip\hfill\qedbox\medbreak\fi\rm}

\def\Smallfonts{}

\def\cite#1{{\bf[#1]}}
\def\Em#1{{\it #1\/}}
\def\Bib#1\par{\bigbreak\bgroup\centerline{#1}\medbreak\parindent30pt
 \parskip2pt\frenchspacing\par}
\def\endBib{\par\egroup}
\newdimen\Overhang
\def\rf#1{\par\noindent\hangafter1\hangindent=\parindent
     \setbox0=\hbox{[#1]}\Overhang\wd0\advance\Overhang.4em\relax
     \ifdim\Overhang>\hangindent\else\Overhang\hangindent\fi
     \hbox to \Overhang{\box0\hss}\ignorespaces}

\def\Coordinates{\bigbreak\bgroup\parindent=0pt\obeylines}
\def\endCoordinates{\egroup}
\let\eval=\Overline
\let\skipaline=\medskip

\H Almost Convex Groups and the Eight Geometries

M. Shapiro and M. Stein\footnote{*}{Both authors thank the NSF for
support.} 

{\Smallfonts \narrower \HHH Abstract

If $M$ is a closed Nil geometry 3-manifold then $\pi_1(M)$ is almost
convex with respect to a fairly simple ``geometric'' generating set.
If $G$ is a central extension or a ${\Bbb Z}$ extension of a word
hyperbolic group, then $G$ is also almost convex with respect to some
generating set.  Combining these with previously known results shows
that if $M$ is a closed 3-manifold with one of Thurston's eight
geometries, $\pi_1(M)$ is almost convex with respect to some
generating set if and only if the geometry in question is not
Sol.\par}

\HH Introduction.

In \cite{C} , Cannon introduced the notion of a finitely generated
group being \Em{almost convex}   with respect to a given generating
set.  This property is  formulated in terms of the geometry
of the Cayley graph, and gives a simple and efficient algorithm for
constructing the Cayley graph.  One class of groups one would like to
study in terms of this property is the class of fundamental groups of
closed 3-manifolds carrying one of Thurston's eight geometries.  (For
an account of these, see \cite{Sc} .)  In most cases, the answer is
already known.  If $G=\pi_1(M)$ where $M$ is a Riemannian
$3$-manifold whose universal cover is $S^3$, ${\Bbb H}^3$, or $S^2
\times {\Bbb R}$, then $G$ is word hyperbolic and thus almost convex
with respect to any finite generating set.  If $M$ is covered by
Euclidean $3$-space, $E^3$, then Cannon \cite{C} shows that $G$ is
almost convex.  Now Cannon {\it et al.}~\cite{CFGT} show that if $M$
is a compact quotient of  Sol, then $G$ is not almost convex with
respect to any generating set.

For the remaining three geometries, it is known that $G$ is almost
convex with respect to some generating set for special cases.
 Specifically, if the
universal cover is Nil and $M$ fibers over a torus, then $G$ is
almost convex with respect to a standard generating set \cite{Sh1}.  
If the universal cover is $\widetilde {PSL_2{\Bbb R}}$ and
$M$ fibers over a closed orientable surface, then $G$ is almost
convex with respect to a standard generating set \cite{Sh2}. 
And if the universal cover is ${\Bbb H}^2 \times {\Bbb R}$ and $G=H
\times {\Bbb Z}$ where $H$ is a hyperbolic surface group, then $G$ is
almost convex with respect to any split generating set.  Note that in
all the remaining cases, $G$ is a finite index supergroup of the
 Nil,  $\widetilde {PSL_2{\Bbb R}}$ and ${\Bbb H}^2 \times {\Bbb R}$ cases
we have just discussed.  Now almost convexity is not a
commensurability invariant.  Indeed, as Thiel has shown \cite{T}, it
is not a group invariant, but does in fact depend on generating set.
Still, it has long seemed likely that the remaining Nil,  $\widetilde
{PSL_2{\Bbb R}}$ and ${\Bbb H}^2 \times {\Bbb R}$ groups  would turn out to
be almost convex.  We will show this to be true, thus proving

\th  Theorem 1

Suppose $G$ is the fundamental group of a closed 3-manifold $M$
carrying one of Thurston's eight geometries.  Then $G$ is almost
convex with respect to some generating set if and only if $M$ is not a
Sol geometry manifold. \endth

This paper is organized as follows.  Section 1 contains background
and definitions.  In Section 2 we establish almost convexity for
the remaining Nil groups.  In Section 3 we establish almost
convexity for the remaining $\widetilde {PSL_2{\Bbb R}}$ and ${\Bbb H}^2
\times {\Bbb R}$ cases.  In the course of this, we prove the following
general result.

\th  Theorem 2

Suppose that $H$ is word hyperbolic, that $A$ is finitely generated
abelian, and that $$1 \longrightarrow A \longrightarrow G \buildrel
\pi \over\longrightarrow H \longrightarrow 1 $$ is a central
extension.  Then there is a finite generating set  $\cal G$ for $G$ so
that $G$ is almost convex with respect to $\cal G$. \endth

\HH 1. Background and definitions.

Let $G$ be a finitely generated group and $\cal G$ a finite set, and
$a \mapsto \eval a$ a map of $\cal G$ to a monoid generating set
$\eval{\cal G} \subset G$.  We use ${\cal G}^*$ to denote the free
monoid on $\cal G$.  We refer to the elements of ${\cal G}^*$ as
\Em{words}.  A word $w=a_1 \ldots a_n$ is said to have \Em{length}
$n$.  This is denoted $\ell(w)=n$.  The map of $\cal G$ into $G$
extends to a unique monoid homomorphism which we denote by $w \mapsto
\eval w$.  We will assume that $\cal G$ is supplied with an involution
$a \mapsto a^{-1}$ and that this involution respects inverses in $G$,
that is, $\eval {a^{-1}}= (\eval a)^{-1}$.  In all cases that we will
consider, the map of $\cal G$ into $G$ is an injection.  Thus we can
make the following convention:  any list of generators will be taken
to include the inverses of those listed.  Thus a definition such as
${\cal X}=\{x,y\}$ will mean ${\cal X}=\{x^{\pm 1},y^{\pm 1}\}$.  If
a generating set contains some element $\rho$ so that $\eval \rho$
has order two, we will take $\rho=\rho^{-1}$.

Given such a group $G$ and $\cal G$, we can form the \Em{Cayley graph}
of $G$ with respect to $\cal G$.  This is a directed labelled graph
$\Gamma = \Gamma_{\cal G}(G)$.  The vertices of $\Gamma$ are the
elements of $G$.  There is a directed edge $(g,a,g')$ from $g$ to $g'$
with label $a$ exactly when $g'=g \eval a$ with $a \in {\cal G}$.
Since $\eval{\cal G}$ generates $G$, $\Gamma$ is connected.  One turns
$\Gamma$ into a metric space by declaring each edge isomorphic with
the unit interval and taking the induced path metric.  We denote this
metric by $d(\cdot , \cdot) = d_{\cal G}(\cdot , \cdot)$.  This, in
turn, gives each  element $x\in \Gamma$ a \Em{length}, $\ell(x)=
\ell_{\cal G}= d_{\cal G}(1,x)$.  As usual, we take the \Em{ball of
radius $r$}, $B(r)$ to be $\{x \in \Gamma \mid \ell(x) \le 1\}$. 

Each edge path in $\Gamma$ is labelled by a unique element of ${\cal
G}^*$.  We identify each word with the edge path it labels starting at
$1 \in G$.  We take this path to be parameterized with unit speed, and
extend each word $w$ to a map of $[0,\infty )$ by setting $w(t)=\eval
w$ for $t \ge \ell(w)$.  The translate of the path $w$ by the group
element $g$ is denoted $gw$.  This is the path based at $g$ bearing
label $w$. We say that $w$ is a \Em{geodesic} if $w|_{[0,\ell(w)]}$ is
an isometry.  Equivalently, $w$ is a geodesic if $\ell(w)=\ell(\eval
w)$.  We say that $w$ is a $(\lambda, \epsilon)$ \Em{quasigeodesic} if
for every subword $u$ of $w$, $\ell(u) \le \lambda \ell(\eval u)+
\epsilon$. 

Following \cite{C}, we say that $G$ is \Em{almost convex} $(m)$ with
respect to $\cal G$ if there is a constant $K(m)$ with the following
property: if $\ell_{\cal G}(g)=\ell_{\cal G}(g')=n$ and $d_{\cal
G}(g,g') \le m$ then there is an edgepath $p$ in $\Gamma$ which runs
from $g$ to $g'$, lies inside $B(n)$ and has length bounded by
$K(n)$.  We say that $G$ is \Em{almost convex} with respect to $\cal
G$ if it is almost convex $(m)$ with respect to $\cal G$ for all $m$.
It is a result of \cite{C} that if $G$ is almost convex $(2)$ with
respect to $\cal G$, then $G$ is almost convex with respect to $\cal
G$.  We will say that $G$ is \Em{almost convex} if there is some
$\cal G$ so that $G$ is almost convex with respect to $\cal G$.  

In Section 3, we will need some standard results concerning word
hyperbolic groups.  For an account of these, see, for example,
\cite{Sho}.  We say that $G$ is \Em{word hyperbolic} if there is a
generating set $\cal G$ and a constant $\delta$ so that if $\alpha$,
$\beta$ and $\gamma$ are geodesic edge paths forming a triangle in
$\Gamma_{\cal G}(G)$, then if $p$ is any point on $\alpha$,
$d(p,\beta\cup\gamma) \le \delta$.  In fact, the existence of such a
$\delta$ is independent of choice of $\cal G$.  Given $G$ and $\cal
G$, we say $D=\{r_1, \ldots , r_k\} \subset {\cal G}^*$ is a
\Em{Dehn's algorithm} if for each $r_i \in D$, we have $\eval {r_i}=1$
and if for any $w \in {\cal G}^*$, if $\eval w = 1$ then  $w$ contains
more than half of some $r_i \in D$ as a subword.  In fact, the
existence of a Dehn's algorithm can be taken as a definition of a word
hyperbolic group.  That is, $G$ is word hyperbolic if and only if for
any generating set $\cal G$, there is a Dehn's algorithm $D \subset
{\cal G}^*$.  Given a Dehn's algorithm $D$, we say a word $w$ is
\Em{$D$-reduced} if $w$ does not contain more than half of any word in
$D$.   By standard methods, one can check that given a Dehn's
algorithm $D$, there are $\lambda$ and $\epsilon$ so that $D$ reduced
words are $(\lambda, \epsilon)$ quasigeodesics.  It is a standard
hyperbolic result that for any distance $m$, there is a constant
$k(m)$ so that if $u,v$ are $(\lambda,\epsilon)$ quasigeodesics with
$d(\eval u, \eval v) \le m$ then each point of $u$ lies within
distance $k(m)$ of $v$ and {\it vice versa}.

\HH 2. Nil manifold groups.

It is shown in \cite {Sh1} \  that if $M$ in a closed 3-manifold with Nil
geometry, and $M$ fibers over a torus then $\pi_1(M)$ is almost
convex.  More specifically, 

\th Theorem 3

Let $$N=\langle x, y \mid [[x,y],x]=[[x,y],y]=1\rangle,$$
$$N^e=\langle x,y,z \mid [x,z]=[y,z]=1,~[x,y]=z^e \rangle.$$ 
We take 
$${\cal N}=\{x,y\} \subset N,$$ 
$${\cal N}^e=\{x,y,z\} \subset N^e.$$  
Then $N$ is almost convex with respect to $\cal N$.  For $e\ge 1$,
$N^e$ is almost convex with respect to ${\cal N}^e$. \endth

\pf Proof

The assertions about $N$ and $N^e$ for $e>1$ are proven in \cite {Sh1} .
The assertion about $N^1$ uses the following lemma.

\th Lemma 4

Suppose $\cal G$ and $\cal H$ are generating sets for $G$ such that
$G$ is almost convex with respect to $\cal G$, and there exists $k$
such that for all $g \in G$, $| \ell_{\cal G}(g) -  \ell_{\cal H}(g)
| \le k$.  Then $G$ is almost convex with respect to $\cal H$. \endth

\pf Proof

Let $\lambda =  \max (\{ \ell_{\cal G}(\eval{h})\mid h\in {\cal
H}\} \cup  \{ \ell_{\cal H}(\eval{g})\mid g\in {\cal G}\}) $.
Suppose $\ell_{\cal H}(h)=\ell_{\cal H}(h')=n$ and $d_{\cal H}(h,h')
\le 2$.  We must find a path of bounded length from $h$ to $h'$ lying
inside the $\cal H$ ball of radius $n$.  We can assume that $n\ge
2k+\lambda+1$, for otherwise, $h$ and $h'$ are connected by a path of
length at most $4k+2\lambda +2$.  Let $a_1 \dots a_n$ and $b_1 \dots
b_n$ be $\cal H$ geodesics for $h$ and $h'$.  Then $g=\eval{a_1\dots
a_{n-2k-\lambda -1}}$ and $g'=\eval{b_1\dots b_{n-2k-\lambda -1}}$ both
lie inside the $\cal G$ ball of radius $n-k-\lambda -1$, and lie
within $\cal G$ distance $\lambda (4k+2\lambda+4)$ of each other.
Since $G$ is almost convex with respect to $\cal G$, they are
connected by a $\cal G$ path $p$ of length at most $K$ which lies
entirely within the $\cal G$ ball of radius $n-k-\lambda -1$.  The
path $p$ is easily turned into an $\cal H$ path $p'$ of $\cal H$
length at most $\lambda K$ whose vertices all lie in the $\cal G$ ball
of radius $n-k-1$.  In particular, $p'$ lies entirely within the $\cal
H$ ball of radius $n$, and so does the path from $h$ to $h'$ labelled
$a_n^{-1}\dots a_{n-2k-\lambda}^{-1} p' b_{n-2k-\lambda}\dots b_n$.
\qed

Continuing with the proof of Theorem 3, we first note that $N$ is
isomorphic to $N^1$.  We now observe that if $g\in N <N^e$ then 
$$\ell_{{\cal N}^e}(g) \le \ell_{\cal N}(g) \le \ell_{{\cal N}^e}(g)
+ 20.$$
The first inequality is obvious.  To see the second, notice that a
geodesic cannot contain both $z$ and $z^{-1}$, for these could be
commuted together and deleted.  Nor can any geodesic contain more than
$25$ $z$'s or $z^{-1}$'s.  For these could be commuted to the end of
the word and replaced with $[x^5,y^5]^{\pm 1}$, thus shortening the
word.  Finally, check that for $-25 \le i \le 25$, $\ell_{\cal N}
(\eval {z^i}) \le 20$.  Now applying Lemma 4 completes the proof
that $N$ is almost convex with respect to ${\cal N}^1$.\qed

We have seen that if $M$ is a Nil manifold that fibers over the torus,
then $\pi_1(M)$ is almost convex with respect to the above generating
sets.  We would like to extend this result to the general case.   Now,
every Nil manifold fibers\footnote{*}{More properly, we should say
``Seifert fibers,'' but we will not keep up this distinction.}  over a
$2$-dimensional Euclidean orbifold.  In fact, the list of orbifolds
$E$ which occur in this role is quite restrictive.  To see this, first
recall from \cite{Sc} that there are no orientation reversing
isometries of Nil.  This implies that $E$ can not have any reflector
curves, for each reflector curve must lift to an orientation reversing
element of $\pi_1(M)$.  Thus the underlying surface of $E$ must be
closed and its singularities (if any) are all cone points.  Now $E$
must have $0$ as its orbifold Euler characteristic.  Thus, if $E$ is
orientable, then $E$ is either the torus or one of $S(2,2,2,2)$,
$S(2,4,4)$, $S(3,3,3)$, $S(2,3,6)$.  (These denote the sphere with
cone points of the orders listed.)  If $E$ is not orientable, then
$E$ is either the Klein bottle or $P(2,2)$, the projecitve plane with
two cone points of order $2$.

In each of these cases, if $G = \pi_1(M)$, we have the diagram

$$\matrix { & & & & &1 & &1 & & \cr
	    & & & & &\downarrow & &\downarrow & & \cr
	    &1 &\rightarrow &{\Bbb Z } &\rightarrow &N^e &\rightarrow &{\Bbb Z}^2 &\rightarrow &1 \cr
	    & & & & &\downarrow & &\downarrow & & \cr
	    &1 &\rightarrow &{\Bbb Z} &\rightarrow &G &\rightarrow &\pi_1^{\rm orb}(E) &\rightarrow &1 \cr
	    & & & & &\downarrow & &\downarrow & & \cr
	    & & & & &Q &\cong &Q & & \cr
	    & & & & &\downarrow & &\downarrow & & \cr
	    & & & & &1 & &1 & & \cr} $$

If $E$ is orientable,  each of the ${\Bbb Z}$ kernels is central.  In
the case where $E$ is not orientable, each element of $\pi_1^{\rm
orb}(E)$ acts trivially or nontrivially on the ${\Bbb Z}$ kernel
depending on whether it is orientation preserving or reversing.

Our strategy for all of these groups is to use the extension 
$$1 \rightarrow N^e \rightarrow G \rightarrow Q \rightarrow 1$$
to lift the almost convexity of $N^e$ up to $G$.  We first establish
the lemmas we will need.

\th Lemma 5

Suppose 
$$1 \rightarrow H \buildrel i \over\rightarrow G \buildrel p
\over\rightarrow Q \rightarrow 1$$
 with $Q$ finite.  Suppose $\cal H$ and $\cal G$ are generating sets
for $H$ and $G$ respectively so that ${\cal H}\subset {\cal G}$, $i$
is an isometry, and the elements of ${\cal G}\setminus  {\cal H}$
permute the elements of $\cal H$ when acting by conjugation.  Then
there is a finite set $T\subset ({\cal G} \setminus
 {\cal H})^*$ so that every element of $G$ has $\cal G$ geodesic
lying in ${\cal H}^*T$. \endth

\pf Proof

We will take $T = \{ t \in ({\cal G}\setminus {\cal H})^* : \ell (t)
\le \# Q \}$.

Let $g \in G$ and suppose $w$ is a geodesic with $\eval w =g$. Since
the elements of ${\cal G}\setminus {\cal H}$ permute the elements of
$\cal H$, it is easy to see that we can find $w'$ so that $\eval w =
\eval {w'}$, $\ell(w) =\ell(w')$ and $w' = h_1 \dots h_m g_1 \dots
g_n$ where $h_i \in \cal H$ for $1 \le i \le m$ and $g_i \in {\cal
G}\setminus {\cal H}$ for $1\le i\le n$.  If $n\le \# Q $, we are
done.  If not, then by the pigeon hole principle, there are $j$ and
$k$, $1 \le j < k \le n$ so that $p(\eval {g_1 \dots g_j}) = p(\eval
{g_1 \dots g_k})$.  Then $\eval {g_{j+1}\dots g_k} \in H$.  Since $i$
is an isometry, we can replace $g_{j+1} \dots g_k$ by an expression of
equal length in ${\cal H}^*$.  Continuing with these two processes
produces a geodesic of the desired form. \qed

\th Corollary 6

Suppose $$ 1 \rightarrow H \buildrel i\over\rightarrow G \buildrel p
\over \rightarrow Q \rightarrow 1$$ with $Q$ finite.  Suppose $\cal H$
and $\cal G$ are generating sets for $H$ and $G$ so that ${\cal
H}\subset {\cal G}$, the elements of ${\cal G}\setminus {\cal H}$
permute the elements of $\cal H$ when acting by conjugation, and $i$
is an isometry.  If $H$ is almost convex with respect to $\cal H$,
then $G$ is almost convex with respect to $\cal G$. \endth

\pf  Proof

Suppose that $g,g'\in G$, $\ell(g)=\ell(g')$ and $d(g,g')\le 2$.  We
invoke the previous lemma and perhaps interchange $g$ and $g'$ to find
geodesics $g=\eval{ut}$, $g'=\eval{u'v't'}$, where $u, u'v' \in {\cal
H}^*$, $t,t' \in T$, and $\ell(v't')=\ell(t)\le \# Q$.  Consequently,
$\eval u, \eval{u'} \in H$, with $d(\eval u,\eval{u'}) \le 2 \# Q +2$.
Using the fact that $H$ is almost convex, there is a path of length at
most $K=K(2\#Q+2)$ connecting $\eval u$ to $\eval {u'}$ inside the
ball of radius $\ell(u)$ in $H$.  This path lies inside the ball of
radius $\ell(u)$ in $G$.  But this gives us a path of length at most
$K+2\# Q$ connecting $g$ to $g'$ inside the ball of radius $\ell(g)$.
\qed

Now we will use  Corollary 6 to show that when $E$ is the Klein
bottle,  $S(2,2,2,2)$, $S(2,2,4)$ or $P(2,2)$, $G$ is almost convex.
First we note that in these cases, the action of $Q$ on ${\Bbb Z}^2$
preserves the standard generating set.  Now in each of these cases, we
can find a generating set for $G$ by lifting the orbifold generating
set for $\pi_1^{\rm orb}(E)$ and appending $z$.  Conjugation by these
generators of $G$ preserves a generating set of $N^e$ of the form
$${\cal N}^e_s=\{x,xz^{i_1}, \ldots , xz^{i_a}, y, yz^{j_1}, \ldots
,yz^{j_b}, z^{k_1}, \ldots ,z^{k_c} \}.$$ We shall call such a
generating set \Em{saturation} of ${\cal N}^e$.  In order to use
Corollary 6, we need to establish

\th Theorem 7

$N^e$ is almost convex with respect to any saturation ${\cal
N}^e_s$ of $\cal N$.\endth 

This follows from Lemma 4 together with  the following 

\th Lemma 8

Let ${\cal N}^e_s$ be a saturation of ${\cal N}^e$.  Then there is a
constant $K$ so that for all $g \in N^e$, $|\ell_{{\cal N}^e_s}(g) -
\ell_{{\cal N}^e}(g) | \le K$. \endth

In order to do this, we need to recall some facts about geodesics in
$N$ in the generating set $\cal N$, established in \cite{Sh1}.  We
review this since we will need to use these methods later for the
$(3,3,3)$ and $(3,6,6)$ groups.  The viewpoint followed in \cite{Sh1}
is to use the short exact sequence $$ 1 \rightarrow {\Bbb Z} \rightarrow
N \buildrel p \over\rightarrow {\Bbb Z}^2 \rightarrow 1 $$ to view words
in ${\cal N}^*$ as lifts of paths in the Cayley graph of ${\Bbb Z}^2$
based at the identity.  Two such paths  $\alpha$ and $\beta$ represent
the same element of $N$ if and only if they end at the same point in
${\Bbb Z}^2$, and the concatenation $\alpha \beta^{-1}$ encloses zero
signed area.  For a given $g=(a,b)\in {\Bbb Z}^2$, we let $B_g(n)=
B_g^{n - \ell(g)}= B(n)\cap p^{-1}(g)$.  (Recall that $B(n)$ is the
ball of radius $n$ in $\Gamma$.  Since $p^{-1}(g)$ consists only of
group elements, so does $B_g(n)= B_g^{n - \ell(g)}$.)

For each $g=(a,b)\in {\Bbb Z}^2$, we identify $p^{-1}(g)$ with ${\Bbb Z}$
by carrying $x^ay^b[x,y]^t$ to $t$.  We observe 
\item{(1)} Geodesics for elements in $p^{-1}(0,0)$ project to closed
paths; in fact, to simple closed curves.  

If you can enclose $n$ squares with a loop of length $l$, you can
enclose $n-1$ squares with a loop of length at most $l$.  This shows
that $B_{(0,0)}^n$ is an interval.  In fact, 
\item {(2)} For any $g \in {\Bbb Z}^2$, $B_g^n$ is an interval.

We will call the elements of $B_g^n$ at the extremes of the
interval \Em{extremals}.  A geodesic or a projection of a geodesic
for an extremal will also be called extremal.  We also observe
\item{(3)} The number of squares enclosed by a closed extremal is
strictly monotone in the length of the extremal.  

>From this we deduce that

\item{(4)} Closed extremals are rectangles.  Indeed, their sides
differ in length by at most one.  We call such a rectangle an
\Em{almost square}.

For suppose $r$ is any simple closed curve.  Let $R$ be the unique
minimal rectangle containing $r$.  $R$ encloses at least as many
squares as little $r$, but $R$ is no longer than $r$.  This is easy
to see:  first notice that minimality of $R$ ensures that $r$ meets
all sides of $R$.  But now given two points on successive sides of
$R$, then the subpath of $R$ that connects them is geodesic, so it is
no longer than the corresponding subpath of $r$.  Now, since the
largest rectangles for a fixed perimeter are almost squares, (4)
follows. 

We are now prepared to describe a sublanguage of the geodesics in $N$
which contains at least one geodesic for each element of $N$ and
suffices for our purposes.  

First notice that for each $g \in {\Bbb Z}^2$, all words representing
elements of $p^{-1}(g)$ differ in length from $g$ by an even amount.
This follows from the fact that all relators are of even length.

\item{(G1)} For each $g \in {\Bbb Z}^2$, the elements of $B^0_g$ are all
represented by geodesics which project to geodesics in ${\Bbb Z}^2$.
(Indeed, every geodesic in ${\Bbb Z}^2$ occurs in this role.)  

\item{(G2)} Any geodesic for an extremal element projects to a subpath
of an almost square.

\item{(G3)} For every element of $B^{2k}_g$, $k \ge 1$, there is a
geodesic whose projection first follows the almost square given by an
extremal of $B^{2k-2}_g$, then crosses via a single edge to follow the
almost square given by an extremal of  $B^{2k}_g$.

We are now ready to prove Lemma 8.

\pf Proof of Lemma 8

Recall our generating sets ${\cal N}^e=\{x,y,z\}$, and
$${\cal N}^e_s=\{x,xz^{i_1}, \ldots , xz^{i_a}, y, yz^{j_1}, \ldots
,yz^{j_b}, z^{k_1}, \ldots ,z^{k_c} \}.$$
Choose $g \in N^e$.  Let $w$ be a geodesic word in ${\cal N}^e_s$ with
$\eval w = g$.  Now let $w'$ be the word in ${{\cal N}^e}^*$ obtained
by deleting all generators of type $z^a$ in $w$ and replacing each
$x^{\pm 1}z^b$ and $y^{\pm 1}z^c$ generator by $x^{\pm 1}$ and $y^{\pm
1}$ respectively.   Let $g'=\eval{w'} \in N <N^e$.  Now, $g=g'z^{\pm
t}$ where $0 \le t \le k \ell(w)$ with
$$k=\max\{|i_1|,\ldots , |i_a|, |j_1|, \ldots , |j_b|, |k_1|, \ldots
|k_c|\}.$$ 
Notice that $[x^m,y^m]=[x,y]^{m^2}=z^{em^2}$. Hence the ${\cal
N}^e$-length of $z^{\pm t}$ is at most $4\Bigl\lceil \sqrt{t/e}\Bigr\rceil +e
\le 4\sqrt{t/e}+e+4$.  So now if $v\in \{x,y\}^*$ is an  $\cal N$
geodesic for $g' \in N$, then 
$$\ell(w) \le \ell(v)+4\sqrt{t/e}+e+4 \le \ell(v)+4 \sqrt{{k\ell(w)
\over e}}+e+4.$$
So if $\ell(w)=\ell_{\cal G}(g)$ is sufficiently large, we will have
$$1/2 \ell(w) \le \ell(v) \le \ell(w).$$ 
Now $g=vz^{\pm t}$ where $t \le k\ell(w) \le 2k \ell(v)$.  But we can
assume that the projection of $v$ looks like one of the projections
described in (G1) --- (G3), and hence has a straight piece of length
$\ell(v)/4$.  So by inserting either $x$ and $x^{-1}$ or $y$ and
$y^{-1}$, around the piece of $v$ corresponding to this side, we
increase the word length by $2$, but increase enclosed area by
$\ell(v)/4$ squares.  Continuing in this way,  and introducing at most
$\lceil 8k/e\rceil$ such pairs of generators and at most $e$ $z$'s, we
produce a word for $g$.  This has increased length by at most $\lceil
16k/e\rceil+e$.  Thus we have
$$\ell_{{\cal N}^e_s}(g) \le \ell_{{\cal N}^e}(g)\le \ell(v)+\lceil
16k/e\rceil +e \le \ell_{{\cal N}^e_s}(g) +\lceil 16k/e\rceil +e.$$
This proves Lemma 8 with $K=\lceil 16k/e \rceil +e$.\qed

We are now ready to establish the almost convexity of the Nil groups
$G=\pi_ 1(M)$ such that $M$ fibers over $E$ where $E$ is the Klein
bottle, $S(2,2,2,2)$, $S(2,4,4)$ or $P(2,2)$.  We will use a
generating set of the form ${\cal G} = {\cal N}^e_s \cup {\cal S}$
where ${\cal S}$ is a lift to $G$ of a standard generating set for
$\pi_1^{\rm orb}(E)$, and ${\cal N}^e_s$ is a saturated generating set
preserved by $\cal S$.   We observed that there are saturated
generating sets preserved by $\cal S$, and that $N^e$ is almost convex
with respect to any saturated generating set.  So all that remains is
to choose a saturated generating set ${\cal N}^e_s$ preserved by $\cal
S$ so that the inclusion $N^e \hookrightarrow G$ is geodesic with
respect to the generating sets ${\cal N}^e_s$ and ${\cal G}={\cal
N}^e_s\cup{\cal S}$. 

We begin by noting that the inclusion ${\Bbb Z}^2\hookrightarrow
\pi_1^{\rm {orb}}(E)$ is geodesic with respect to the generating sets
${\cal X}=\{x,y\}$ and ${\cal E}=  \{\rho,x,y\}$, $\{a,b,c,d,x,y\}$,
$\{p,q,r,x,y\}$, $\{a,b,\rho,x,y\}$  depending on whether $E$ is the
Klein bottle,  $S(2,2,2,2)$ or $S(2,4,4)$ or $P(2,2)$ (In each case,
$\rho$ is orientation reversing, $a,b,c,d,p$ are rotations of order
$2$ and $q$ and $r$ are rotations of order $4$.)  We
spell out the argument for $E= S(2,2,2,2)$.  We need to show that no
geodesic in $\{x,y\}^*$ may be shortened by rewriting it in the
$\{a,b,c,d,x,y\}$ generating set.  So suppose $w'$ is an $\{x,y\}$
geodesic which can be shortened to a $\{a,b,c,d,x,y\}$ geodesic $w$.
We can push the $\{a,b,c,d\}$ letters all to the right without
changing the number of $\{x,y\}$ letters. Thus we assume can assume
$w=uv$, where $u$ is composed of $\{x,y\}$ letters and $v$ is composed
of $\{a,b,c,d\}$ letters.  But now, since $Q={\Bbb Z}_2$, any substring
$v$ of length $3$ or more includes a substring which evaluates into
${\Bbb Z}^2$.  One checks that words in $\{a,b,c,d\}^*$ of length less
than $3$ which evaluate into ${\Bbb Z}^2$ are not shorter than the
corresponding $\{x,y\}$ words for the elements they represent.  By
this process we can eliminate all $\{a,b,c,d\}$ letters from $w$
without increasing length.  Since $w'$ was assumed to be an $\cal X$
geodesic, this contradicts the assumption that $w$ is shorter than
$w'$.  

The argument for the other cases proceeds similarly.  When $E$ is the
Klein bottle, we check $\{\rho\}$ words of length $2$.  Since
$\rho^2=x$ the result is immediate.  For $E=S(2,4,4)$, one checks all
$\{p,q,r\}$ words of length at most $4$.  Similarly for the case
$E=P(2,2)$, we must check $\{a,b,\rho\}$ words of lengths at most $4$.

We now wish to lift this to $N^e < G$.  Choose a saturated generating
${\cal N}'_s$ for $N^e$ which is preserved by conjugation by $\cal S$.
Suppose the inclusion $N^e \hookrightarrow G$ is not geodesic with
respect to ${\cal N}'_s$ and ${\cal G}'={\cal N}'_s \cup {\cal S}$.
As above any failure of this inclusion to be geodesic can be observed
in some short word $v \in {\cal S}^*$ with $\ell(v) \le \#Q$.
Consider the projection of $v$ into $\pi_1^{\rm{orb}}(E)$.  Since the
inclusion down below is geodesic, there is an $\{x,y\}$ word $v'$
which evaluates to this projection and $\ell(v') \le \ell(v)$.  We now
consider $v'$ as a word in the generators for $N^e$.  Note that $\eval
v = \eval {v'z^k}$ for some $k$.  Further, only finitely many such $k$
occur since there are only finitely many short $v$'s.  Let $K$ be the
maximum of the absolute values of those $k$ which arise.  We take
$${\cal N}^e_s = \{az^t \mid a \in {\cal N}'_s, |t| \le K\}.$$ Clearly
the inclusion of $N^e$ into $G$ is geodesic with respect to ${\cal
N}^e_s$ and ${\cal G}={\cal N}^e_s\cup{\cal S}$, and ${\cal N}^e_s$ is
preserved by the action of $\cal S$.

Thus we have shown that if $G=\pi_1(M)$ is a Nil group so that $M$
fibers over the Klein bottle, $S(2,2,2,2)$, $S(2,4,4)$ or $P(2,2)$,
there is a generating set $\cal G$ so that $G$ is almost convex with
respect to $\cal G$.

It is now easy to see that $G$ is almost convex with respect to the
generating set ${\cal N}^e\cup{\cal S}$ or  ${\cal N}\cup{\cal S}$ if
this latter generates $G$.  This is because every element of $G$ has a
$\cal G$ geodesic in which there are at most $\#Q$ $\cal S$ letters
and these occur at the end.  But now Lemma 8 tells us that replacing
the ${\cal N}^e_s$ part by the corresponding ${\cal N}^e$ or $\cal N$
geodesic increases length by at most a bounded amount. Thus by Lemma
4, $G$ is almost convex with respect to this reduced generating set.

We are now prepared to pursue the case where $M$ fibers over
$S(3,3,3)$ or $S(2,3,6)$.  We wish to pursue the same program.
However, the action of our finite quotient no longer preserves the $x$
and $y$ directions in the plane, and hence does not preserve the
generating set $\{x,y\}$.  Rather, it preserves hexagonal symmetry,
and hence preserves generating sets of the form $\{x,y,xy\}$.  (We may
think of our $x$ and $y$ directions as labelling nonadjacent rays in
this hexagonal symmetry.)  Thus we will need to work with generating
sets of the form $${\cal X} =\{x,y,t\} \subset {\Bbb Z}^2,$$ $${\cal N}
=\{x,y,t\} \subset N,$$ $${\cal N}^e =\{x,y,t,z\} \subset N^e.$$ (In
each case, we take $t=xy$.)  To carry out our program, we must first
prove

\th Theorem 9

$N$ is almost convex with respect to ${\cal N}$. \endth

\th Theorem 10

$N^e$ is almost convex with respect to  ${\cal N}^e$.\endth

Any generating set of the form 
$${\cal N}^e_s= \{x,xz^{i_1}, \ldots xz^{i_a},y,yz^{j_1}, \ldots
yz^{j_b}, t,tz^{k_1}, \ldots tz^{k_c}, z,z^{l_1},\ldots ,
z^{l_d}\}$$ 
is called a \Em{saturation} of ${\cal N}^e$.

\th Corollary 11

$N^e$ is almost convex with respect to any saturation ${\cal
N}^e_s$.\endth 

Once again, we study $N$ by studying the projection of paths in $N$
into ${\Bbb Z}^2$.  This time, however, we have taken the generating set
${\cal N}=\{x,y,t\}$ and its projection $\cal X$ in ${\Bbb Z}^2$.
 The Cayley graph of ${\Bbb Z}^2$ with respect to $\cal X$ is the
$1$-skeleton of the tessellation of the plane by equilateral
triangles.  Now $\pi_1^{\rm orb}S(3,3,3)$ acts as the orientation
preserving subgroup of the symmetries of this tessellation.  There are
two orbits of triangles under this action, those with boundary label
$xyt^{-1}$ (which we color white) and those with boundary label
$tx^{-1}y^{-1}$ (which we color black).  A circuit around a white
triangle lifts to the identity in $N$. A circuit around a black
triangle lifts to $\eval {[x,y]}$ in $N$.  Consequently, a closed
curve in the plane lifts to $\eval{[x,y]^n}$, where $n$ is the signed
number of black triangles enclosed.  As before, we can use this
technology to treat $N$ as equivalence classes of based edgepaths in
the Cayley graph of ${\Bbb Z}^2$.  Then after replacing the words
``almost square'' with ``almost regular hexagon'' (1) --- (4) and (G1)
--- (G3) hold.  (A hexagon is an \Em{almost regular hexagon} if no two
of its sides differ in length by more than $2$.)  

>From this description of geodesics in $N$ we can deduce

\th Lemma 12

Suppose that $g,g' \in N$ and that $\ell(g)=\ell(g')=n$ and that
 $g'=g \eval r$ with $\ell(r) \le 2$.  Then one of the following
occurs 
\item {1)} Both $g$ and $g'$ are small, i.e., each has length
less than $50$.  
\item {2)} There are standard geodesics for $g$ and $g'$ whose
projections in ${\Bbb Z}^2$ lie in an $25$ neighborhood of each other
and $50$-fellow travel.  
\item {3)} The projection of $gr$ crosses an axis. If it crosses
exactly one axis, then there is $g''$ so that $\ell(g'')=n$, $d(g,g'')
\le 2$, the projection of $g''$ lies on the axis, $d(g',g'') \le 2$,
and 2) above holds for each of the pairs $(g,g'')$ and $(g'',g')$.
\endth

The above estimates are not sharp.

\pf Sketch of the proof of Lemma 12

We suppose that $g$, $g'$, and $r$ are given as above, and suppose
that $\alpha'$ and $\beta'$ are geodesics for $g$ and $g'$.  If $g\in
B^0$, then the projection of $\alpha'$ lies entirely in the
parallelogram determined by the tessellation of the plane with corners
at $(0,0)$ and $p(g)$.  In this case we can further demand that the
projection of $\alpha'$ have a specific form: we demand that it stay
along one side of the parallelogram as long as possible.  In this case
we take $\alpha$ to be the path in the plane given by the projection
of $\alpha'$.  If $g \notin B^0$, then $\alpha'$ can be taken to lie
(in projection) close to the boundary of an almost regular hexagon. We
take $\alpha$ to be the path along the boundary of this hexagon
determined by $\alpha'$.  We choose $\beta$ similarly.  

One can then perform a systematic enumeration of the possibilities
using the following facts: 
\item {(E1)} Each of $\alpha$ and $\beta$ is either a geodesic of the
special form we have given or consists of between one and six sides
of an almost regular hexagon.
\item{(E2)} Each of $\alpha$ and $\beta$ starts along one of the $6$
axial directions.
\item{(E3)} If $\alpha$ or $\beta$ lies along an almost regular
hexagon, it turns either clockwise or counterclockwise.

Without loss of generality, one picks an axial direction for $\alpha$
and quickly eliminates the cases in which $\alpha$ and $\beta$ turn
in opposite directions.

One has the following data:
\item {(D1)} The lengths of $\alpha$ and $\beta$ differ by at most $4$.
\item{(D2)} The end points of $\alpha$ and $\beta$ are separated by a
path $s$ of length at most $4$.
\item{(D3)} The path given by $\alpha'r\beta'^{-1}$ in projection
encloses $0$ signed area, where area is measured by the number of
black triangles enclosed.  
\item{(D4)} The projection of $\alpha'$ and $\alpha$  are either
identical or lie separated by a straight strip of width $1$.  In the
case where $\alpha$ is more than $2$ sides of an almost regular
hexagon, we can take $\alpha$ to have the same endpoint as the
projection of $\alpha'$, and similarly for $\beta$.  Thus, when 
$\alpha$ and $\beta$ are more than $2$ sides of an almost regular
hexagon, $s$ is simply the projection of $gr$.

One then proceeds to examine each of these cases using elementary
Euclidean geometry.  In each of the cases where at least one of
$\alpha$ or $\beta$ does not consist of $6$ sides of an almost
regular hexagon, we find that we are in case 1) or 2) of the Lemma.
That is, either $g$ and $g'$ are both short, or there are geodesic
paths for $g$ and $g'$ which fellow travel in projection.

The case where both $\alpha$ and $\beta$ are $6$ sides of an almost
regular hexagon is more interesting.  If $\alpha$ and $\beta$ both end
in the interior of a common sextant (i.e., in the region between two
axial rays), then $\alpha$ and $\beta$ must start in the same axial
direction and have approximately the same size.  Thus, they fellow
travel.  On the other hand, if, say, $\alpha$ ends in the interior of
a sextant and $\beta$ ends on an axis which defines that sextant,
$\alpha $ and $\beta$ need not fellow travel in projection.  For here,
$\alpha$ may start out in (say) the $t$ direction, while $\beta$ can
start out in the $y$ direction, so that we have, say,
$\beta=y^it^jx^jy^{-j}t^{-j}x^{-j}$, with $i<j$.  However, the fact
that $\beta$ starts in the $y$ direction and ends on the $y$ axis
allows us to ``slide'' $\beta$ in the following manner: for any
$\delta \le i$, the path
$y^{i-\delta}t^jx^jy^{-j}t^{-j}x^{-j}y^\delta$ evaluates to the same
element of $N$ as $\beta$.  By taking $\delta=j$ we replace $\beta$
with an equivalent geodesic which starts in the same direction as
$\alpha$.  Their almost regular hexagons are approximately the same
size, and thus once again,  $\alpha$ and $\beta$ fellow travel.

Finally, if the path $gs$ from the end of $\alpha$ to the end of
$\beta$ crosses exactly one axis, we need only note that there is an
element $g''$ which is within distance $2$ of each of $g$ and $g'$,
has the same length as these and projects to a point lying on this
axis. \qed

\pf Proof of Theorem 9 from  Lemma 12

Let $\ell(g)=\ell(g')=n$ and $d(g,g')\le 2$, say $g'=g\eval r$, with
$\ell(r) \le 2$.   Thus, we are in the situation of Lemma 12.  We must
exhibit a path of bounded length which connects $g$ to $g'$ and lies
in $B(n)$.

If situation 1) of the Lemma holds, then we can connect $g$ to $g'$
by going through the identity, and this path has length at most $100$.

Suppose situation 2) of the Lemma holds.  We shall also assume that
$g$ and $g'$ have length at least $900$ for otherwise, we proceed as
above.  (We remind the reader that our estimates are in fact very
crude!).  Let $w=uv$ and $w'=u'v'$ be geodesics for $g$ and $g'$ with
$ \ell(v)=\ell(v')=900$.  Let $q$ be the lift to $N$ of a path
connecting $p(\eval u)$ to $p(\eval{u'})$. We can assume $\ell(q) \le
50$.  Then the path given by $v^{-1}qv'r^{-1}$ in projection encloses
an area $A$ with $|A|\le (50)(900)$.  We let $\gamma$ be a geodesic
for $z^{-A}$.  This has length at most $4 \sqrt{A}+4<900$.  We then let
$P$ be the path based at $g$ bearing the label $v^{-1}q\gamma v'$.
It is easy to see that this connects $g$ to $g'$ and lies inside
$B(n)$. 

We now suppose that we are in situation 3) of the Lemma.  If the
projection of $gr$ crosses only one axis, then we can perform the
previous process twice, once connecting $g$ to $g''$, and once
connecting $g''$ to $g'$.  But in fact, $r$ has length at most $2$, so
the projection of $gr$ can cross at most $2$ axes, and this happens
only when $g$ and $g'$ project to points within distance $1$ of
$(0,0)$.  In this case we repeat the previous method twice and are
done. \qed

\pf Proof of Theorem 10

This follows from Theorem 9 by observing that as usual an ${\cal N}^e$
geodesic in $N^e$ can have very few $z$'s.  Then if $g,g' \in N^e$
with $\ell(g)=\ell(g')$ and $d(g,g')\le 2$, we write geodesics
$w=uz^m$ and $w'=u'z^{m'}$ for $g$ and $g'$, and $|m|$ and $|m'|$ are
bounded, and $u$ and $u'$ are free of $z$'s.  Without loss of
generality we can assume $|m|\le |m'|$.  We let $r$ be the terminal
segment of $u$ of length $|m'|-|m|$, so that $u =u''r$ with
$\ell(u'')=\ell(u')$.  Now the almost convexity of $N$ gives us a
path $Q$ of bounded length connecting $\eval{u''}$ to $\eval{u'}$ inside
the ball of radius $\ell(g)-|m'|$.  The path we seek is the one which
starts at $g$ and is labelled $z^{-m}r^{-1}Qz^{m'}$. \qed

The Corollary now follows by the methods of Lemma 8.

We now have all the tools in place that we used for the proof in the
square case.  The proof in the present case follows along exactly the
same lines.  One checks that the embedding of ${\Bbb Z}^2$ into the
appropriate orbifold groups is geodesic with respect to the generating
set ${\cal X}=\{x,y,t\}$ for ${\Bbb Z}^2$ and the appropriate orbifold
generating sets ${\cal S}\cup{\cal X}$.  The action of $\cal S$
preserves $\cal X$.  Once again, this lifts to the inclusion $N^e
\hookrightarrow G$ giving almost convexity with respect to a
generating set of the form ${\cal N}^e_s \cup{\cal S}$, where ${\cal
N}^e_s$ is a saturated generating set.  As before, this gives almost
convexity with respect to the generating set ${\cal N}^e \cup {\cal
S}$ or ${\cal N} \cup {\cal S}$ if this latter generates.

\HH 3. Central extensions of word hyperbolic groups.

We will now prove Theorem 2.

\pf Proof of Theorem 2

We take $\rho : H \times H \rightarrow A$ to be the cocycle defining
$G$.  Thus we can identify $G$ with the set $A \times H $ endowed with
the multiplication $(a,h)(a',h')= (a+a'+\rho(h,h'), hh')$.  Let $\cal
G'$ be a generating set for $G$.  Then ${\cal H}=\pi {\cal G'}$ is a
generating set for $H$.  We will extend $\cal G'$ to a generating set
$\cal G$ so that $G$ is almost convex with respect to $\cal G$.

\th  Lemma  13

Suppose $\cal G'$ is a generating set for $G$ with $\pi \cal G' = \cal
H$.  Suppose that $D$ is a Dehn's algorithm for $H$ with respect to
$\cal H$.   Then there are sets ${\cal G''}\subset \pi^{-1}{\cal H}$
and  ${\cal A'} \subset A$ so that if ${\cal A'} \subset {\cal A}$ and
${\cal G} = {\cal G'}\cup {\cal G''}  \cup {\cal A}$, then for any $g
\in G$ there is a geodesic $w \in {\cal G}^*$ so that $\eval w = g$
and $\pi w$ is $D$-reduced. \endth

\pf Proof

We start by taking 
$${\cal G}''=\{(1,h) \mid h \in {\cal H}\}.$$  

We enlarge $D$ if necessary to make sure it is closed under inversion
and cyclic permutation.  We wish to lift the words in $D$ to $A$.  For
each $d = h_1 \dots h_k \in D$, we take $\tilde d = (1,h_1) \dots
(1,h_k) \in A$.  For each $d= h_1 \dots h_k \in D$, let 
$$R_d = \{ (a_1,h_1)\dots (a_k,h_k) : (a_i,h_i)\in {\cal G}'\cup{\cal
G''}~{\rm{for}}~1\le i \le k\}.$$
For each $r=(a_1,h_1)\dots (a_k,h_k) \in R_d$ and each $i \le k$, let
$a_i(r)=a_1+ \dots + a_i$.  We take 
$${\cal A'} = \{ a_i(r)+\tilde d : r\in R_d, \ d \in D \}.$$

We take an arbitrary geodesic $w''=(a_1,h_1)\dots(a_p,h_p)$ for
$g \in G$.  Since $A$ is central in $G$, we may replace $w''$ with $w'$ of
the form 
$$w'=(a_1,1)\dots (a_i,1)(a_{i+1}, h_{i+1})\dots (a_p,h_p)$$
where $\eval {h_j} \neq 1$ for $j > i$.  We then have $\pi (w')=
h_{i+1} \dots h_p$.  If this is $D$-reduced we are done.  If not,
there are $m$ and $n$ with $i<m<n \le p$ so that $h_m \dots h_n$ is
more than half a relator in $D$.  Hence $\eval{h_m \dots h_n} = \eval
{h'_m \dots h'_{n'}}$ with $n'<n$ and with $d=h_m \dots h_n (h'_m
\dots h'_{n'})^{-1} \in D$.  We then have 

$$\eqalign{\eval{(a_m,h_m)\dots (a_n,h_n)} &= \eval {(a_m+\dots +
a_n,1)(1,h_m)\dots (1,h_n)}\cr 
&=  \eval {(a_m+\dots + a_n+\tilde d,1)(1,h'_m)\dots (1,h'_{n'})}.\cr}$$
By construction, $(a_m+\dots + a_n+\tilde d,1) \in {\cal G}$, so this
last expression lies in ${\cal G}^*$, and since $n'<n$, it is no
longer than the first expression.  Thus, we may use it to replace the
first expression in $w'$.  This reduces the length of $\pi w'$, so 
continuing in this way produces a geodesic $w$ whose projection is
$D$-reduced.\qed

Notice that this means that $\cal A'$ is a generating set for $A$.
For suppose $w$ is as guaranteed by Lemma 13.  If $\eval w \in A$,
then $\eval {\pi w}=1 \in H$, and a $D$-reduced path for the identity
is the empty path.  In particular $A$ is geodesic in $G$, that is,
given a generating set $\cal G$ of the form ${\cal G}'\cup{\cal
G}''\cup{\cal A}$ as above, the inclusion of $A$ into $G$ is an
isometry.

We continue with our proof that $G$ is almost convex.  We must choose
our generating set.  Recall that $D$-reduced words are
$(\lambda,\epsilon)$ quasigeodesics in $H$, and that there is a $k$
so that if two $(\lambda,\epsilon)$ quasigeodesics end at most
distance $2$ apart, then they lie in $k$ neighborhoods of each other.
We will take 
$${\cal A} = {\cal A'} \cup \{ a \in A : a = \eval r {\rm{~with~}} r
\in ({\cal G'}\cup{\cal G''})^* {\rm{~and~}} \ell (r) \le  \lambda
(2k+3) + \epsilon +2k+3 \}.$$ 
We take ${\cal G} = {\cal G'}\cup {\cal G''}  \cup {\cal A}$.

We now suppose that $g, g' \in G$ with $\ell(g)=\ell(g')=n$ and $g'=g
\eval q $ with $\ell(q) \le 2$.  We write $g = \eval w$ where $w = uv$
so that each letter of $u$ projects to $1 \in H$ and no letter of $v$
does.  Similarly, we write $g'=\eval {w'}$ with $w'=u'v'$.  Notice
that $d(\eval{\pi v},\eval{\pi v'}) \le 2$.  Since $\pi v$ and $\pi
v'$ are both $D$-reduced, they are $(\lambda, \epsilon)$
quasigeodesics lying in $k$ neighborhoods of each other.

Suppose first that $\ell(v) \ge k+1$, say $v=xy$ with $\ell(y)=k+1$.
Then $\ell(\eval {ux})=n-k-1$ and $d(\eval{\pi x },\pi w') \le k$.
Choose $z$ so that $\ell(z)\le k$ and $\pi z$ labels a path from
$\eval {\pi x}$ to a point on $\pi v'$.  We will suppose this point to
be $\eval {\pi x'}$ where $v'= x'y'$.  Notice that the path $uxz$
stays inside the ball of radius $n-1$ Since $\pi v'$ is a $(\lambda,
\epsilon)$-quasigeodesic, it follows that $\ell(y')\le \lambda (2k+3)
+ \epsilon $.  Hence the path $y'qy^{-1}z$ has length at most
$\ell(y')\le \lambda (2k+3) + \epsilon +2k+3$.  It projects to a
closed path, so, in particular $\eval {y'qy^{-1}z} = (a,1) \in {\cal
A}$.  We now have $g \eval {y^{-1}z(a,1)} = \eval {u'x'}$.  But the
path labelled $y^{-1}x(a,1)$ based at $g$ stays inside the ball of
radius $n$.  Clearly the path labelled $y'$ based at $\eval{u'x'}$
also stays inside this ball.  Thus, the path labelled
$y^{-1}z(a,1)^{-1}y'$ runs from $g$ to $g'$ staying inside the ball of
radius $n$.  Its length is clearly bounded.

We must now check that we can produce such a path when $\ell(v) \le
k$.  Suppose $\ell(v) \ge 1$.  Then ``backing up along $\pi v$'' takes
us to $1 \in H$, and we can perform the same argument as above taking
$y=v$ and $z$ to be trivial.

We are now left with the case where $\ell(v)=0$ and by symmetry, we
may assume $\ell(v')=0$.  In that case $g$ and $g'$ lie in the abelian
group $A$.  We have seen that $A$ is geodesic in $G$, so we are done.
\qed

This gives the following

\th Scholium 14

Suppose that $H$ is word hyperbolic and that

$$1 \longrightarrow {\Bbb Z}  \longrightarrow G  \longrightarrow H
\longrightarrow 1.$$
Then there is a generating set $\cal G$ so that $G$ is almost convex
with respect to $\cal G$. \endth

\pf Proof

There are only two actions on ${\Bbb Z}$, namely the trivial action and
the action which inverts elements of ${\Bbb Z}$.  Thus at the possible
cost of inverting elements of ${\Bbb Z}$, we can move each of these to
the beginning of a word.  Now this process cannot increase length.
Consequently, each element of $G$ has a geodesic in which every $\Bbb
Z$ generator appears the beginning.  Now one can proceed as above.

\th Corollary 15

Let $M$ be a closed 3-manifold with $\widetilde{\mathop{\rm PSL}_2\Bbb
R}$ or ${\Bbb H}^2 \times {\Bbb R}$ geometry.  Then there is a generating
set $\cal G$ so that $\pi_1(M)$ is almost convex with respect to $\cal
G$.  \endth

\pf Proof

In this case 
$$1 \longrightarrow {\Bbb Z} \longrightarrow \pi_1(M)\longrightarrow H
\longrightarrow 1, $$ 
where $H$ is the orbifold fundamental group of a hyperbolic surface
orbifold.  Since this is necessarily word hyperbolic, the result
follows. \qed

This completes the proof of Theorem 1.

\Bib{References}
\MaxReferenceTag{CFGT}

\rf{C} 
J.~Cannon, Almost convex groups, Geometrae  Dedicata {\bf 22}, 197---210
(1987).

\rf{CFGT} J.~Cannon, W.~Floyd, M.~Grayson, W.~Thurston,  Solvgroups
are not almost convex, Geometrae Dedicata {\bf 31}, 291---300 (1989).

\rf{Sc}
P.~Scott, The Geometries of three-manifolds, Bull.\  London Math.\
Soc.\ {\bf 15}, 401---487 (1983).

\rf{Sh1} M.~Shapiro, A Geometric approach to the almost convexity and
growth of some nilpotent groups, Mathematische Annalen, {\bf 285},
601---624 (1989).

\rf{Sh2} M.~Shapiro, Growth of a $\widetilde{\mathop{\rm PSL}_2\Bbb
R}$ manifold group, to appear in Mathematische Nachrichten.

\rf{Sho}  H.~Short, ed., Notes on word hyperbolic groups, in Group
Theory from a Geometric Viewpoint, E.~Ghys, A.~Haefliger,
A.~Verjovsky eds., World Scientific.

\rf{T} C.~Thiel, Zur Fast-Konvexit\"at einiger nilpotenter Gruppen,
Bonner Mathematische Schriften, 1992.

\endBib

\Coordinates

\noindent City College\\

\noindent New York, NY 10031
\skipaline
\noindent Ohio State University\\

\noindent Columbus, OH 43210
\endCoordinates

\bye